\documentclass[11pt]{article}

\usepackage[utf8]{inputenc}
\usepackage{titling}
\usepackage{lipsum}
\usepackage[space]{grffile}
\usepackage{graphicx}
\usepackage{url}
\usepackage{amsmath}
\usepackage{amssymb}
\usepackage{float}

\usepackage{verbatim}
\usepackage{xcolor}
\usepackage{pbox}

\begin{document}

\vspace*{-2cm}
\noindent \texttt{\small {\em The Mathematical Gazette}, (Cambridge University Press), 107:29, 2023.}
\vspace*{1cm}
\begin{center}
{\Large\bf A Cautionary Example Relating to the\\ 
               Interpretation of Numerical Results}\\
               \vspace{6pt}
               {\large Jeffrey Uhlmann}\\
               Dept.\ of Electrical Engineering and Computer Science\\
               University of Missouri - Columbia
               
\end{center}

\section{Introduction}

The availability of software tools \cite{wolfram,mathematica,octave,maple,matlab,sage} and multi-function calculators has unquestionably had an impact on how students interpret numerical solutions to mathematical problems. However, these tools can also lead to an overly-casual attitude about how to interpret the effects of numerical precision. For example, students quickly learn to interpret a numerical result of $0.99999999$ as being exactly $1$, or $3.54e$-$16$ as being zero. The problem, of course, is that they may become so habituated to disregarding low-order terms that they fail to recognize results that are {\em nearly} an integer -- {\em but are not}. The following trig expression (in radians) provides a cautionary example: 
\begin{equation}
    \arcsin(1+\sin(11))-\sin(11) 
    \label{udiff}
\end{equation}

\begin{figure*}
	\centering
		\includegraphics[scale=1.0]{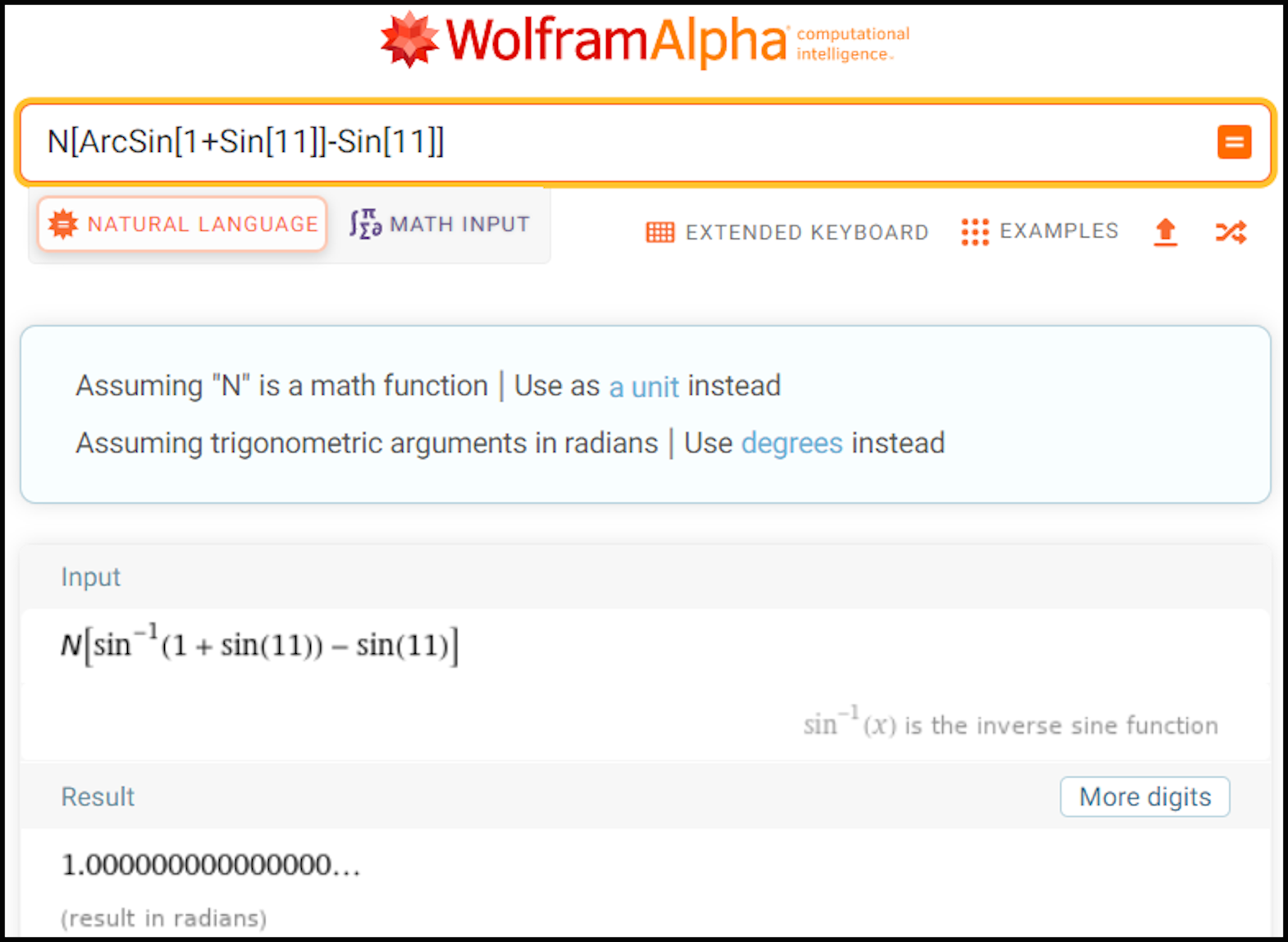}
	\caption{\footnotesize Result of the evaluation of Equation \ref{udiff} using {\em WolframAlpha} \cite{wolfram} with default precision (and defined in radians, not degrees). Note that the ellipsis (...) can reasonably be taken to imply that the sequence of decimal zeros continues indefinitely, or it can be taken to imply there are additional significant digits (possibly nonzero) that are not shown.}
        \label{wolfram1}
\end{figure*}

\begin{figure*}
	\centering
		\includegraphics[scale=1.6]{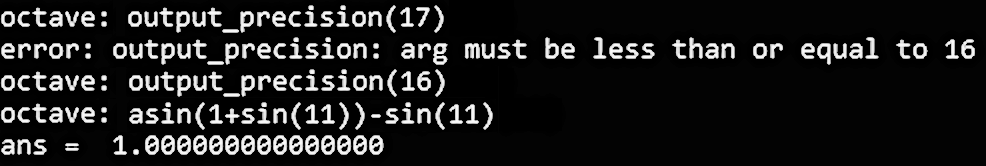}
	\caption{\footnotesize Result of the evaluation of Equation \ref{udiff} using {\em Octave} \cite{octave} with maximum output/display precision.}
	\label{octave}
\end{figure*}

\begin{figure*}
	\centering
		\includegraphics[scale=1.2]{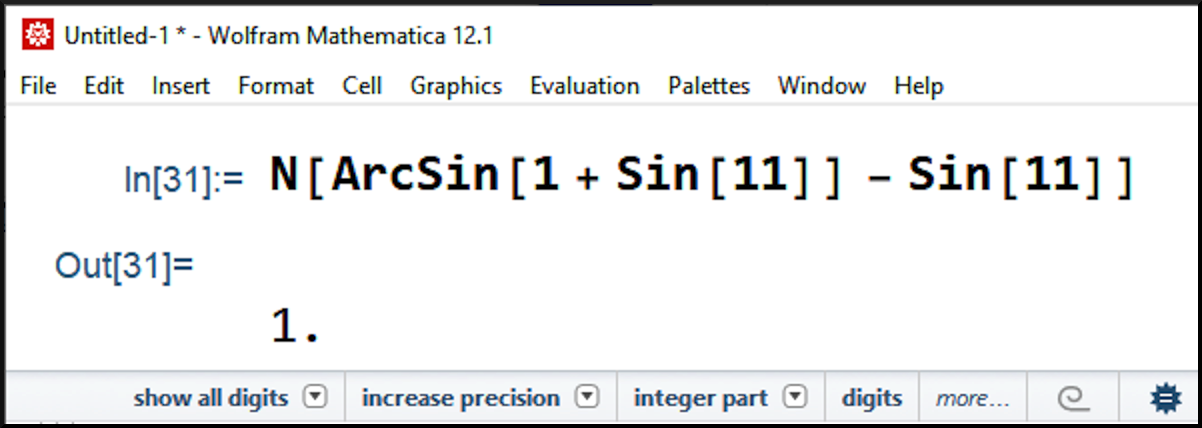}
	\caption{\footnotesize Result of the evaluation of Equation \ref{udiff} using {\em Mathematica} \cite{mathematica} with default precision. It may not be clear from the form of the result ({\bf 1.})  whether the answer is precisely the integer $1$ or has been rounded.}
	\label{mathematica}
\end{figure*}

Figure \ref{wolfram1} shows the result for Equation \ref{udiff} using {\em WolframAlpha}, which is probably the most widely used online mathematical tool. As can be seen, it seems to indicate that the expression evaluates to $1$. Figures \ref{octave} and \ref{mathematica} show results obtained from {\em Octave/Matlab} and {\em Mathematica}, respectively, that are similarly $1$ to within the precision of the display. In fact, a student could be forgiven for concluding from the {\em Mathematica} result in Figure \ref{mathematica} that the answer given is {\em exactly} $1$, i.e., not just {\em nearly} $1$ up to some level of precison.

\begin{figure*}
	\centering
		\includegraphics[scale=1.1]{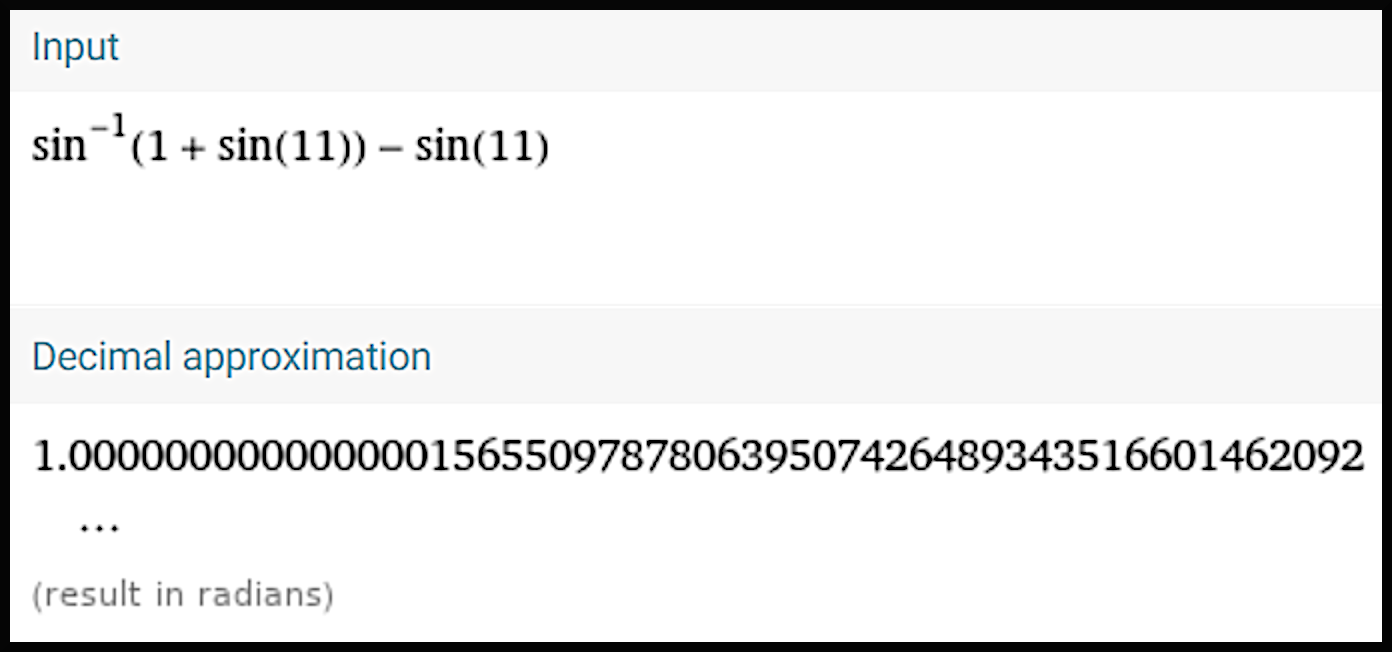}
	\caption{\footnotesize Result of the evaluation of Equation \ref{udiff} using {\em WolframAlpha} with increased displayed precision. The shown digits are entirely accurate, i.e., the expression does not equal 1, but at first glance the natural inclination might be to assume that the nonzero decimal digits are simply due to  {\em numerical inaccuracy}.}
       \label{wolframz}
\end{figure*}

If a student were to increase the displayed precision using {\em WolframAlpha} or {\em Mathematica}, the result would be what is shown in Figure \ref{wolframz} with zeros out to 16 decimal places. A student might very well conclude that digits beyond that long sequence of zeros are just due to numerical inaccuracy. This is a natural assumption to make based on practical experience with calculators and other common mathematical tools. However, students are likely to be very surprised when they discover that the number shown in Figure \ref{wolframz} is in fact accurate to the precision shown, i.e., that the expression of Equation \ref{udiff} does {\em not} equal $1$.

\section{Considerations on Numerology}

There are multiple factors that may contribute to students drawing a mistaken inference from the example of the previous section. First is the fact that some calculators display spurious degrees of precision. This can occur when equations are solved to a precision less than the default display precision. One might expect a few spurious low-order digits to be of little consequence -- let alone be noticed. However, student homework problems are typically designed to be amenable to solving by hand and often have integer solutions.  Thus, students become habituated to seeing slight deviations from known integer results due to numerical imprecision from their calculators. It is no wonder, then, that such experience could lead to false presumptions when real-world problems {\em happen} to produce near-integer solutions. 

A natural response to this concern might be to point out that such events are relatively rare. After all, the chances of a solution to a `random' problem being within 6 or 7 digits of an integer is on the order of one in a million. Yes, that's roughly correct, but what if a solution {\em happens} to be near $\pi$, $2\pi$, $\pi/2$, $\sqrt{2}$, or some other commonly-encountered special number? Is it likely that a student might assume a result beginning with `3.14159...' is identically $\pi$ even though discrepancies appear at the 6th decimal place? Possibly so. Thus, opportunities for numerological misadventure are somewhat greater than might be estimated at first glance. 

The risk of numerological seduction is further complicated by the fact that proximity to an integer, or some other special number, could potentially hint that `something deeper' may be at play. For example, in 1859 Charles Hermite recognized that $\exp(\pi\sqrt{163})$ (now commonly referred to as Ramanujan's constant \cite{barrow}) is a near-integer to 12 decimal places, and this has inspired subsequent efforts to explain {\em why} this is the case. Yes, there could be something deeper, but if its near-integer status is purely coincidental, then search for an `explanation' can only yield another numerological chimera that is no more revealing than the existence of a line that connects two given points.

\section{Discussion}

One way to reduce susceptibility to mistakes of the kind discussed in this article is to avoid using calculators that display spurious digits of precision. Educators can also emphasize to students that many modern tools (e.g., {\em WolframAlpha / Mathematica}) never display spurious digits, so if a displayed result deviates from an integer or other special number, then it truly is {\em not} equal to that special number. More generally, a full discussion of this topic with students is likely to be stimulating and provoke greater sensitivity to the possiblity that tantalizing features of a given solution may be purely accidental.

\vspace{1cm}

\noindent {\footnotesize {\em Acknowledgement}: The author would like to thank Warren D.\ Smith for an interesting conversation about ``specialness'' of the Hermite-Ramanujan number.}

\end{document}